\documentclass[12pt]{amsart}
\usepackage{somedefs,tracefnt,latexsym,,rawfonts,epsfig,amsfonts,amssymb,latexsym,enumerate,amscd}
\def\cal{\mathcal}
\def\Bbb{\mathbb}
\def\g{\gamma}
\def\G{\Gamma}
\def\r{\rangle}
\def\l{\langle}
\def\t{\times}
\def\p{\partial}
\input psfig
\newtheorem{prop}{Proposition}[section]
\newtheorem{thm}{Theorem}[section]
\newtheorem{exm}{Example}[section]
\newtheorem*{thm1}{Main Theorem}
\newtheorem*{thm2}{Reduction Theorem}
\newtheorem*{thm3}{The Loop Theorem}
\newtheorem{lemma}{Lemma}[section]
\newtheorem{cor}{Corollary}[section]

\newtheorem{defn}{Definition}[section]

\newtheorem{rem}{Remark}[section]
\numberwithin{equation}{section}
\begin{document}
\title[The Farrell-Jones isomorphism conjecture] {The Farrell-Jones
isomorphism conjecture for $3$-manifold groups}
\author[S.K. Roushon]{S.K. Roushon}
\address{School of Mathematics\\
Tata Institute\\
Homi Bhabha Road\\
Mumbai 400005, India}
\email{roushon@math.tifr.res.in}
\urladdr{http://www.math.tifr.res.in/$\sim$roushon/paper.html}
\date{March 13, 2006}
\begin{abstract} We show that the Fibered Isomorphism Conjecture
(FIC) of
Farrell and Jones corresponding to the stable topological pseudoisotopy
functor is true for the fundamental groups of a large class of
$3$-manifolds. We
also prove that if the FIC is true for irreducible $3$-manifold
groups then it is true for all $3$-manifold groups. In fact, this follows
from a more general result we prove here, namely we show that if the FIC
is true for each vertex group of a graph of groups with trivial edge
groups then the FIC is true for the fundamental group of the graph of
groups. This result is part of a program to prove FIC for the 
fundamental group of a graph of groups where all the vertex and edge
groups satisfy FIC. A consequence of the first result gives a partial
solution to a problem in the problem list of R. Kirby. We also deduce that
the FIC is true for a class of virtually $PD_3$-groups.

Another main aspect of this article is to prove the FIC for all Haken
$3$-manifold groups assuming that the FIC is true for $B$-groups. By
definition a $B$-group contains a finite index subgroup isomorphic to 
the fundamental group of a compact irreducible
$3$-manifold with incompressible nonempty boundary so that each boundary  
component is of
genus $\geq 2$. We also prove the FIC for a large
class of $B$-groups and moreover, using a recent result
of L.E. Jones we show that the surjective part of the FIC is true for any
$B$-group.\end{abstract} 

\keywords{$3$-manifold, virtual fibering, fibered isomorphism conjecture,
Whitehead group, stable topological pseudo-isotopy functor}

\subjclass[2000]{Primary: 57N37, 19J10. Secondary: 19D35.}

\maketitle

\tableofcontents

\section{Introduction}
Algebraic $K$-theory of $3$-manifold groups was considered by Waldhausen
(\cite{W}) leading to some very general
results in algebraic $K$-theory which has far wide applications. Namely,
the results of Waldhausen regarding algebraic $K$-theory of generalized
free products and $HNN$-extensions. The motivation for these celebrated
results was computing Whitehead groups of knots groups and of $3$-manifold
groups in general. This led him to prove the very important result that
the Whitehead group, reduced projective class group and the negative
$K$-groups all vanish for the fundamental group of any Haken $3$-manifold.
Another outstanding result in this direction is due to Farrell and Jones.
They proved that the same result is true for the fundamental group of a
hyperbolic manifold (\cite{FJ1}, \cite{FJ2}) and more generally for
nonpositively curved manifolds (\cite{FJ3}) (of any dimension). Thurston's
Geometrization
conjecture says that an irreducible aspherical $3$-manifold is either
Seifert fibered, hyperbolic or Haken. Plotnick (\cite{Pl}) proved the
vanishing of the above $K$-theoretic groups for non-Haken Seifert
fibered spaces. Thus conjecturally this completes the picture for
computing low dimensional algebraic $K$-theory of aspherical $3$-manifold
groups. In the setting of assembly map these $K$-theoretic vanishing
results are implied by the isomorphism of some assembly map for
pseudoisotopy functor. The more recent Fibered Isomorphism Conjecture
(FIC) of Farrell and Jones predicts that this isomorphism to hold for all
groups (\cite{FJ}). 

In this paper we are concerned about proving the FIC for several classes
of $3$-manifold groups and $PD_3$ groups. At first we establish a general
result showing that the FIC is true for a free product if it is true for
each free summand. From this result and some standard facts in 
$3$-manifold topology we deduce that it is enough to consider irreducible
$3$-manifold groups only. Another aspect of this paper is to show that if
the FIC is true for groups (we call such a group a $B$-group) containing a
finite index subgroup
isomorphic to the fundamental group of compact irreducible
$3$-manifold with incompressible nonempty boundary so that each
boundary component is a surface of genus $\geq 2$, then it is true for
all $3$-manifold groups modulo Thurston's Geometrization conjecture. We
also prove the FIC for a subclass
of $B$-groups. We further show that any torsion free 
$B$-group can be obtained by applying generalized free products or
$HNN$-extensions (amalgamated along infinite cyclic subgroup) on
finitely many members of the above subclass. 

In this paper by the FIC we always mean the Fibered Isomorphism Conjecture 
of Farrell and Jones (\cite{FJ}) corresponding to the stable topological
pseudoisotopy functor. In short it says that computing the
pseudoisotopy functor of $K(\pi,1)$ is equivalent to the
computation of the pseudoisotopy functor of $K(V, 1)$ where $V$ varies
over all virtually cyclic subgroups of $\pi$. In recent times the FIC has
been checked for
several classes of groups and it has become an important fundamental tool
for
computing $K$-theory of groups (see for example  \cite{BFJP}, \cite{FJ},
\cite{FL}, \cite{FR}, \cite{R2}, \cite{R3} etc.). 

In the paper \cite{Q1} F. Quinn mentioned that the technique we use to
prove
FIC in the pseudoisotopy functor case can also be used in the case of
$K$-theory developed in \cite{Q1} to deduce a similar conclusion for
higher $K$-theory with arbitrary coefficient ring.

\section{Statement of the FIC and related results}\label{appendix}

In this section we recall the Fibered Isomorphism Conjecture of 
Farrell and Jones made in \cite{FJ}. The following formulation of the
conjecture is taken from \cite{FL}. 

Let $\cal S$ denotes one of the three functors from the category of
topological spaces to the category of spectra: (a) the stable topological
pseudo-isotopy functor ${\cal P}()$; (b) the algebraic $K$-theory functor 
${\cal K}()$; (c) and the $L$-theory functor ${\cal L}^{-\infty}()$. 

Let $\cal M$ be the category of continuous surjective maps. The objects of 
$\cal M$ are continuous surjective maps $p:E\to B$ 
between
topological spaces $E$ and $B$. And a morphism between two maps $p:E_1\to
B_1$ and $q:E_2\to B_2$ is a pair of continuous maps $f:E_1\to E_2$,
$g:B_1\to B_2$ such that the following diagram commutes. 

There is a functor defined by Quinn \cite{Q} from $\cal M$ to the category
of $\Omega$-spectra which associates to the map $p:E\to B$ the spectrum
${\mathbb H}(B, {\cal S}(p))$ with the property that ${\mathbb H}(B, {\cal
S}(p))={\cal S}(E)$ when $B$ is a single point. For an explanation of
${\mathbb H}(B, {\cal S}(p))$ see [\cite{FJ}, section 1.4]. Also the map
${\mathbb H}(B, {\cal S}(p))\to {\cal S}(E)$ induced by the morphism:
id$:E\to E$; $B\to *$ in the category $\cal M$ is called the Quinn
assembly map. 

$$\begin{CD}
E_1 @>f>> E_2\\
@VVpV    @VVqV\\
B_1 @>g>> B_2
\end{CD}$$

Let $\Gamma$ be a discrete group and $\cal E$ be a $\Gamma$ space which
is universal for the class of all virtually cyclic subgroups of $\Gamma$ 
and denote ${\cal E}/\Gamma$ by $\cal B$. For definition of universal 
space see [\cite{FJ}, appendix]. Let $X$ be a space on which $\Gamma$ acts
freely and properly discontinuously and $p:X\t_{\Gamma} {\cal E}\to
{\cal E}/{\Gamma}={\cal B}$ be the map induced by the projection onto the 
second factor of $X\t {\cal E}$. 

The Fibered Isomorphism Conjecture states that the map $${\mathbb
H}({\cal B}, {\cal S}(p))\to {\cal S}(X\t_{\Gamma} {\cal E})={\cal
S}(X/\Gamma)$$ is a (weak) equivalence of spectra. The equality is induced
from the map $X\t_{\Gamma}{\cal E}\to X/\Gamma$ and using 
the fact that $\cal S$ is homotopy invariant.

Let $Y$ be a connected $CW$-complex and $\Gamma=\pi_1(Y)$. Let $X$ be the
universal cover $\tilde Y$ of $Y$ and the action of $\Gamma$ on $X$ is
the action by group of covering transformation. If we take an aspherical
$CW$-complex $Y'$ with $\Gamma=\pi_1(Y')$ and $X$ is the universal cover
$\tilde Y'$ of $Y'$ then by [\cite{FJ}, corollary 2.2.1] if the FIC is
true for the space $\tilde Y'$ then it is true for $\tilde Y$ also. Thus
throughout the paper whenever we say that the FIC is true for a discrete
group
$\Gamma$ or for the fundamental group $\pi_1(X)$ of a space $X$ we would
mean it is true for the Eilenberg-MacLane space $K(\Gamma, 1)$ or
$K(\pi_1(X), 1)$ and the functor ${\cal S}()$.

Now we recall the relation between the $\Omega$-spectra ${\cal P}()$ and
lower algebraic $K$-theory. This result is proved in \cite{AH}.
$$\pi_j({\cal P}(M))=\begin{cases}
K_{j+2}({\mathbb Z}\pi_1(M))& \text{if}\ j\leq -3,\\
\tilde K_0({\mathbb Z}\pi_1(M))& \text{if}\ j=-2,\\
Wh({\mathbb Z}\pi_1(M))& \text{if}\ j=-1.\end{cases}$$

The following results are used frequently in this paper.

Recall that a wreath product $A\wr B$ of two groups $A$ and $B$ is by 
definition the semi-direct
product $A^B\rtimes B$ with respect to the regular action of $B$ on
$B$. Also by the algebraic lemma of \cite{FR} if $A$ is a normal subgroup 
of $B$ with quotient $C$ then $B$ can be embedded in the wreath product 
$A\wr C$.

\begin{defn} A group $G$ is said to satisfy {\it FICwF} if the FIC is true
for $G\wr K$ for any finite group $K$.\end{defn}

\begin{lemma} \label{lemma A} ([\cite{FJ}, theorem A.8]) If the FICwF is
true for a discrete group $\Gamma$ then it is true for any subgroup of
$\Gamma$.\end{lemma}

\begin{lemma}\label{lemma C} ([\cite{FJ}, proposition 2.2]) Let $f:G\to H$
be a surjective homomorphism. Assume that the FIC is true for $H$ and for
$f^{-1}(C)$ for all virtually cyclic subgroup $C$ of $H$ (including
$C=1$). Then the FIC is true for $G$.\end{lemma}

Lemma \ref{lemma C} also holds if we replace FIC by FICwF.

\begin{prop}\label{prop2.4} [\cite{FJ}, proposition 2.4] The FICwF is true
for virtually poly-$\Bbb Z$ groups.\end{prop}

\section{Statement of Reduction Theorem: the FIC for free products and
reduction to irreducible case}\label{freeproduct}

Let us first recall some basic definitions from $3$-manifold topology. By
a {\it $3$-manifold group} we mean the fundamental group of a
$3$-manifold. A $3$-manifold is called {\it irreducible} if every
embedding of the $2$-sphere in the $3$-manifold extends to an embedding of
the $3$-disc. And a $3$-manifold is called {\it prime} if whenever it is a
connected sum of two $3$-manifolds then one of the summand is homeomorphic
to the $3$-spheres. The only $3$-manifolds which are prime but not
irreducible are the ${\mathbb S}^2$ bundles over the circle ([\cite{He},
lemma 3.13]). And by Kneser's prime decomposition theorem any compact
$3$-manifold can be written as a connected sum of finitely many prime
$3$-manifolds ([\cite{He}, theorem 3.15]). Thus we have the following
lemma.

\begin{lemma}\label{prime} A compact $3$-manifold group is a free product
of finitely many irreducible $3$-manifold groups and a free group of
finite rank.\end{lemma} 

The first theorem of this article is about the FIC for arbitrary groups. 

\begin{thm2}\label{reduction} Let $G_1$ and $G_2$ be two countable groups
and assume
that the FIC is true for both $G_1$ and $G_2$. Then the FIC is true for
the free product $G_1*G_2$.\end{thm2}

\begin{thm} \label{refree} Assume that the FICwF is true for $G_i$
for $i=1,2$. Then the FICwF is true for $G_1*G_2$.\end{thm}

\begin{cor} \label{reductioncor} If the FICwF is true for compact 
irreducible aspherical $3$-manifold groups then it is true for any
$3$-manifold group.\end{cor}

Because of Corollary \ref{reductioncor} and the Remark below in the next 
sections we always consider orientable irreducible aspherical 
$3$-manifolds. 

\begin{rem} {\rm As any nonorientable $3$-manifold is covered by
orientable ones and we always prove the FICwF for $3$-manifold groups,
the nonorientable case follows
using [\cite{FR}, algebraic lemma].}\end{rem} 

\medskip

A more interesting corollary to the Reduction Theorem is the following
application to graph of groups. Recall that a {\it graph of groups}
$\cal G$ is a graph, that is a connected one dimensional $CW$-complex, 
with the following additional data.

\begin{itemize}
\item for each vertex $v$ there is an associated group $G_v$.
\item for each edge $e$ there is an associated group $G_e$ and for the
end
vertices $v^e_i$ and $v^e_f$ of $e$ there are injective homomorphisms
$G_e\to G_{v^e_i}$ and $G_e\to G_{v^e_f}$. 
\end{itemize}

Note that $v^e_i$ and $v^e_f$ could be the same vertex.

In the above definition $G_v$ is called a {\it vertex group} and $G_e$ is
called an {\it edge group}.

For example, a generalized free product $G_1*_HG_2$ can be associated 
with the graph $\cal G$ with two distinct vertices $v_1$ and $v_2$ and an
edge $e$ joining them. The vertex groups are $G_{v_1}=G_1$ and
$G_{v_2}=G_2$. The edge group is $G_e=H$ and the two injections $G_e\to
G_{v_1}$ and $G_e\to G_{v_2}$ are the injections defining the generalized
free product $G_1*_HG_2$. In this case $G_1*_HG_2$ is called the
fundamental group of the graph $\cal G$. Similarly if $G*_H$ is an
$HNN$-extension then the corresponding graph consists of a single vertex
$v$ and an edge $e$ whose both end points are the same vertex. The
following Figure describe the situation in detail.

Let $\cal G$ be a finite graph of groups and $N_{\cal G}=$(number of
vertices of $\cal G$ + number of edges of $\cal G$). Now we proceed to
define the {\it fundamental group} $\pi_1({\cal G})$ of the finite graph
of groups $\cal G$. At first when $N_{\cal G}=1$ then the graph has one
vertex $v$ and no edge. Define $\pi_1({\cal G})=G_v$. By induction
on $N_{\cal G}$ assume that we have defined $\pi_1({\cal G})$ for all
graph $\cal G$ for which $N_{\cal G}\leq n-1$. So let $\cal G$ be a finite
graph and $N_{\cal G}=n$. There are now two cases to consider. 

\vskip 0.5cm

\centerline{\psfig{figure=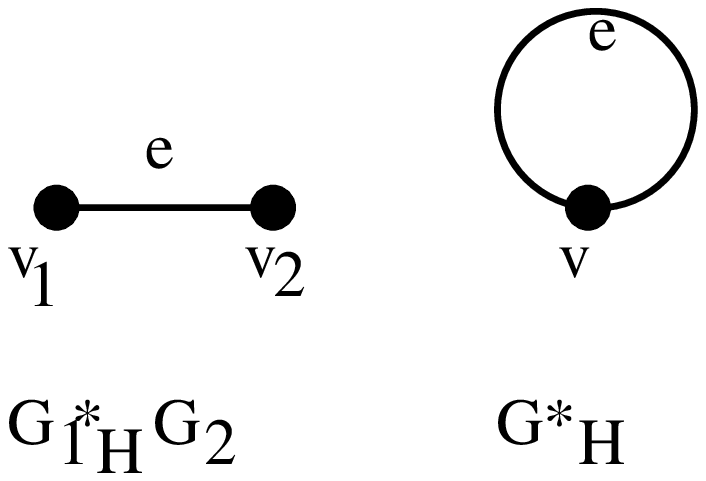,height=3cm,width=5cm}}
%\vskip 0.5cm

\centerline{Figure: graph of groups corresponding to generalized}
\centerline{free product and $HNN$-extension.}
\vskip 0.7cm

\noindent
{\bf Case A.} Let $v$ be a vertex of $\cal G$ with only one edge $e$ 
emanating form it. Let $G_v$ be the group associated to the vertex $v$.
Let ${\cal G}_1$ be the graph obtained by deleting $v$ and $e$ from $\cal
G$. Clearly ${\cal G}_1$ has $N_{{\cal G}_1}=n-2$. Hence by the induction
hypothesis $\pi_1({\cal G}_1)$ is defined. Now define $\pi_1({\cal G})=
\pi_1({\cal G}_1)*_{G_e}G_v$.

\noindent
{\bf Case B.} All the vertices have valency $\geq 2$, that is, there are
at least two edges emanating from every vertex of the graph. Remove one
edge $e$ from the graph and let ${\cal G}_1$ be the resulting graph. Hence
${\cal G}_1$ has $N_{{\cal G}_1}=n-1$ and by induction $\pi_1({\cal G}_1)$
is defined. Then define $\pi_1({\cal G})=\pi_1({\cal G}_1)*_{G_e}$. 

To define the fundamental group of an infinite graph of groups $\cal G$ 
note that $\cal G$ can be written as an increasing union of finite
subgraphs ${\cal G}_i$. Define $\pi_1({\cal
G})=\lim_{i\to\infty}\pi_1({\cal G}_i)$. 

\begin{rem}{\rm Note that in the above definition of fundamental group of
a graph of groups there are many choices involved. But this definition
does our purpose as any two choices define isomorphic groups.}\end{rem}

The Reduction Theorem is for the simplest
nontrivial graph, namely the graph with two vertices and one edge. The
edge group is trivial and the vertex groups are $G_1$ and $G_2$. 

\begin{cor}\label{graphfic} Let $\cal G$ be a graph of groups with
trivial edge groups. Also assume that the graph has countable number of
vertices and edges and each vertex group is countable. If the FICwF is
true for all vertex groups then FICwF 
is true for $\pi_1({\cal G})$.\end{cor}

\section{Statements of theorems and proofs of basic results: the FIC for
$3$-manifold groups}\label{ficmanifold}

To state the theorems we need to recall some more preliminaries from
$3$-manifold topology. An embedded closed $2$-manifold $F$ in a
$3$-manifold is called $2$-{\it sided} if the normal bundle of $F$ in $M$ 
is
trivial, that is homeomorphic to $F\t {\Bbb R}$.  For example any 
orientable $2$-manifold embedded in an orientable $3$-manifold is  
$2$-sided. By a {\it Haken}
$3$-manifold we mean it is irreducible and contains an embedded $2$-sided
$\pi_1$-injective closed $2$-manifold of infinite
fundamental group. Such an embedded $2$-manifold is called an {\it
incompressible} surface. Also we recall that there is a unique
decomposition of any compact Haken $3$-manifold by cutting the
$3$-manifold along finitely many incompressible tori so that each piece is
either Seifert fibered or hyperbolic. This is known as JSJT (Jaco-Shalen,
Johannson and Thurston) decomposition. Recall that a {\it graph manifold}
is either a Seifert fibered space or is a Haken $3$-manifold which
has only Seifert fibered pieces in
its JSJT decomposition or equivalently which is obtained by gluing along
tori boundary components of finitely many Seifert fibered $3$-manifolds. 

Let us first state a theorem which is a consequence of some fundamental 
results by Farrell and Jones for manifolds of any dimension. 

\begin{thm} (\cite{FJ}, \cite{FJ3}) \label{Farrell-Jones} Let $M$ be a 
closed
nonpositively curved Riemannian manifold or a compact surface (may be 
with nonempty boundary). Then the FICwF is true for
$\pi_1(M)$.\end{thm}

\begin{proof} In the closed manifold case the theorem is obtained by 
combining theorem 
0.4 of \cite{FJ3} and the remark following it together with proposition 
2.3 of \cite{FJ}.
In the closed surface case we also have to 
use Proposition \ref{prop2.4}. When 
the surface has nonempty boundary take the double $N$ of $M$. 
Then $\pi_1(M)\wr K$ is a subgroup of $\pi_1(N)\wr K$ for any finite group 
$K$. Now using Lemma \ref{lemma A} and the previous case we complete the 
proof.\end{proof}

The following corollary to the above Theorem is crucial for this paper.

\begin{cor} \label{jsjt} Let $M$ be a closed Haken $3$-manifold such that
there is a hyperbolic piece in the JSJT decomposition of the manifold 
or let $M$ be a compact irreducible $3$-manifold with nonempty
incompressible boundary and there is at least one torus boundary
component. Then the FICwF is true for $\pi_1(M)$.\end{cor}

\begin{proof} In the closed case by \cite{Le} $M$ supports a nonpositively
curved Riemannian metric. Theorem \ref{Farrell-Jones} proves the corollary
in this case.

Secondly let $M$ be compact with incompressible boundary and there is at 
least one torus boundary. Let $M'$ be the double of $M$ along
boundary components of genus $\geq 2$. Then $M'$ is irreducible and with 
incompressible tori boundary components. Using \cite{Le} we deduce that
the interior of $M'$ supports a complete nonpositively curved
Riemannian metric so that near the boundary components the metric is a
product. Now let $M''$ be the double of $M'$. Then $M''$ supports a
nonpositively curved Riemannian metric. Hence by the closed case the FICwF
is true for $\pi_1(M'')$. On the other hand we have the inclusions
$\pi_1(M)<\pi_1(M')<\pi_1(M'')$ which induce $\pi_1(M)\wr G<\pi_1(M')\wr
G<\pi_1(M'')\wr G$. Here $G$ is a finite group. Hence the FICwF is true 
for $\pi_1(M)$ by Lemma
\ref{lemma A}.\end{proof}

Another result we will be using often is the following theorem of Farrell
and Linnell.

\begin{thm} ([\cite{FL}, theorem 7.1]) \label{Farrell-Linnell} Let $I$ be
a directed set, and let $\Gamma_n$, $n\in I$ be a directed system of
groups with $\Gamma=\lim_{n\in I}$; i.e., $\Gamma$ is the direct limit of
the groups $\Gamma_n$. If each group $\Gamma_n$ satisfies the FIC, then
$\Gamma$ also satisfies the FIC.\end{thm}

Let $M$ be a complete Riemannian manifold. Then $M$ is called {\it
A-regular} if the following is satisfied. $$|\nabla^i(R)|\leq A_i$$ where
$A_i$ is a sequence of nonnegative integers and $\nabla^i(R)$ denotes the
$i$-th covariant derivative of the curvature tensor $R$ of the Riemannian
manifold $M$. For example; any compact or homogeneous or locally symmetric
Riemannian manifold is $A$-regular. Some more examples, relevant for this
article, of noncompact $A$-regular Riemannian manifolds are described in
Section \ref{arrm}. 

\begin{defn} {\rm A group $\Gamma$ is called an {\it A-group} if there
is a complete simply connected $A$-regular Riemannian manifold $M$ of
nonpositive sectional curvature and a properly discontinuous action of
$\Gamma$ on $M$ by isometries so that the image of $\Gamma$ in Iso$(M)$ is
virtually torsion free.}\end{defn}

L.E. Jones recently proved that [\cite{Jo}, theorem 1.6] the assembly map
in the statement of the Fibered isomorphism Conjecture induces a
surjective homomorphism on the homotopy group level for torsion
free $A$-groups. Note that the FIC states that this homomorphism should be
an isomorphism. Note that $A$-regularity is a consequence if the 
action in the above Definition is also cocompact. Under this cocompactness
assumption the FIC is proved in \cite{FJ} for any $A$-groups, even if it
has torsion (that is the Theorem \ref{Farrell-Jones}). 

\medskip
\noindent
{\bf Notation.}
Let us denote by $\cal C$ the class of compact irreducible $3$-manifolds
with nonempty incompressible boundary so that each boundary component is a
surface of genus $\geq 2$. 

\begin{defn} \label{bgroupd}{\rm A group $\Gamma$ is called a {\it
B-group} if it contains a finite index subgroup isomorphic to the
fundamental group of a member of $\cal C$.}\end{defn}

In Proposition \ref{agroup} we show that if $M\in {\cal C}$ then
$\pi_1(M)\wr G$ is an $A$-group for any finite group $G$. 

Though it is not yet proved that the FIC is true for all $A$-groups, we
have some partial results in this direction. 

Let $M$ be a compact $3$-manifold with nonempty boundary. A compact
surface $F\subset M$ is called {\it properly embedded} if $\p F\subset \p
M$ and $F-\p F\subset M- \p M$. 

\begin{defn} \label{essential} {\rm A properly embedded annulus $A$ which
is the image  of an embedding $g:(F, \p F)\to (M, \p M)$ is called {\it
essential} if the followings hold.

\begin{itemize}

\item the embedding $g$ is not isotopic (relative to boundary) to an
embedding $f:(F, \p F)\to (M, \p M)$ so that $f(F)\subset \p M$. 

\item $g_*:\pi_1(F)\to \pi_1(M)$ is injective.

\end{itemize}}
\end{defn}

\begin{rem}{\rm The above definition of essential embedding is little
different from standard definition. Usually the second condition is not
assumed in standard terminology and in the first condition it is only
demanded that $g$ is homotopic to a map $f:(F, \p F)\to (M, \p M)$ so that
$f(F)\subset \p M$.}\end{rem}

In Section \ref{bgroup} we prove some results describing properties of
members of $\cal C$. We show in Proposition \ref{building} that, for any 
$M\in {\cal C}$, $\pi_1(M)$ is isomorphic to
the fundamental group of a graph of groups with infinite cyclic edge
groups and whose vertex groups are  
fundamental groups of members of $\cal C$ which contain no essential
annulus.

The following theorem is a first step towards proving the FIC for  
$B$-groups.

\begin{thm} \label{annulus} Let $M\in {\cal C}$. Assume that there is no
essential annulus embedded in $M$. Then the FICwF is true for
$\pi_1(M)$.\end{thm}

The following Proposition describes the main set of examples of $A$-groups
used in this article. We prove this Proposition in Section \ref{arrm}.  

\begin{prop}\label{agroup} Let $M\in {\cal C}$. Then the followings are 
true. 

\begin{itemize}

\item In the JSJT decomposition of $M$ the pieces which contain a
boundary component of $M$ supports hyperbolic metric in the interior.

\item The interior of $M$ supports a nonpositively curved Riemannian 
metric so that near each boundary components the metric has constant $-1$ 
sectional curvature. Also the group $\pi_1(M)\wr G$ is an $A$-group for 
any finite group $G$.
\end{itemize}
\end{prop}

A second important step to prove the FIC for $B$-groups is the following.

\begin{thm} \label{geodesic} Let $M\in {\cal C}$. Consider the Riemannian
metric in the interior of $M$ as given by Proposition \ref{agroup}. 
Assume
that with respect to the hyperbolic metric near each boundary component
the boundary components are totally geodesic. Then the FICwF is true for
$\pi_1(M)$.\end{thm} 

In Examples \ref{examp} we show that there are members of $\cal C$ which
do
not satisfy the hypothesis of Theorems \ref{annulus} and \ref{geodesic}.
We do not know if these are the only such examples.

We are now ready to state the main result of this article.

\begin{thm1} Assume that the FIC is true $B$-groups\footnote{Recently in
\cite{R3} we have shown that the FIC is true for $B$-groups.}. Let $M$ be
a
closed $3$-manifold. Let $H$ be a homomorphic image of
$\pi_1(M)$ satisfying the following properties.

\begin{itemize}
\item $H$ has a finite index nontrivial torsion free subgroup.
\item the FICwF is true for $H$.
\item any infinite cyclic subgroup of $H$ has infinite index in $H$.
\end{itemize}

Then the FICwF is true for $\pi_1(M)$.\end{thm1}

There are several important consequences of the Main Theorem and its
proof. Before we state the consequences we recall some background.

Recall that by Thurston's hyperbolization theorem if a Haken $3$-manifold
contains no incompressible torus then it is hyperbolic and hence in this
case Theorem \ref{Farrell-Jones} applies. Thus in the class of Haken
$3$-manifolds, because of Corollary \ref{jsjt}, the FIC remains to be
proven for closed graph manifolds only.  

The Geometrization conjecture implies that an irreducible
$3$-manifold either has finite fundamental group or it is one of the
following types: Seifert fibered, hyperbolic or Haken.  

By a {\it virtually fibered} $3$-manifold we mean that it has a finite
sheeted cover which fibers over the circle.

It is known that there are aspherical $3$-manifold which are not virtually
fibered. But a well-known conjecture says that every irreducible
$3$-manifold with infinite fundamental group (and hence aspherical) has a
finite sheeted cover which is Haken (in this case the manifold is called 
{\it virtually Haken}). This is called the {\it virtual Haken
conjecture}. 
 
\begin{thm} \label{corollary} Assume that the FIC is true for $B$-groups.  
Then the FICwF is true for $\Gamma$ where $\Gamma$ is isomorphic to the
fundamental group of one of the following manifolds.

\begin{itemize} 
\item compact $3$-manifold with nonempty boundary or
a noncompact $3$-manifold.

\item $3$-manifold which has a finite sheeted cover 
with first Betti number $\geq 2$.

\item nontrivial graph manifold. 

\item virtually fibered $3$-manifold. 

\item virtually Haken $3$-manifold.

\end{itemize}
\end{thm}

Here we should point out that the different items in the above
theorem are not mutually
exclusive. But we state them separately for the importance and
popularity of the individual classes of manifolds and groups.

In this connection we recall that in \cite{R2} we proved the following.

\begin{prop}\label{special}[\cite{R2}, proposition 7.1] Let $M$ be a
compact
$3$-manifold and there is a finite sheeted covering of $M$ which fibers
over the circle with special monodromy
diffeomorphism. Then the FICwF is true for $\pi_1(M)$.\end{prop}

For definition of special diffeomorphism of a surface see [\cite{R2},
definition, p. 4].

Now we state some Theorems where we do not need the assumption that `the
FIC is true for $B$-groups' always. Important instances are Theorems
\ref{seifert} and \ref{graph}. 

\begin{thm} \label{seifert} Let $S$ be a Seifert fibered space.
Then the FICwF is true for $\pi_1(S)$.\end{thm}

An important corollary of Corollary \ref{reductioncor}, Theorem
\ref{corollary} and Theorem
\ref{seifert} is the following.

\begin{cor} Assume that the FIC is true for $B$-groups. Then the FICwF is
true for any $3$-manifold (irreducible or not) group if Thurston's
Geometrization conjecture is
true.\end{cor}

Below, in Theorem \ref{graph} we prove the FIC for some large classes of 
$3$-manifolds and for a class of $PD_3$ groups in Theorem \ref{pdgroup}. 

Note that we do not assume that $M$ is a Haken $3$-manifold in these 
theorems.

Recall that if $M$ is a compact $3$-manifold and there is a surjective
homomorphism from $\pi_1(M)$ to $\mathbb Z$ with finitely generated kernel
then by Stallings' fibration theorem $M$ fibers over the circle and hence
by Theorem \ref{corollary} the FIC is true for $\pi_1(M)$ provided the
FIC is true for any $B$-group. In fact here we deduce a more general
statement showing that the fourth item of Theorem \ref{corollary} is
true for some general class of $3$-dimensional Poincar\'{e} duality
($PD_3$) groups. This theorem is an application of Theorem  
\ref{corollary} and [\cite{Hi}, theorem 1.20] and the well-known result
that a $PD_2$ is a surface group (\cite{EL} and \cite{EM}).

\begin{thm} \label{pdgroup} Assume that the FIC is true for any $B$-group.
Let $\cal G$ be a group and $H$ a finite index, finitely presented
subgroup of $\cal G$. Assume that $H$ is a $PD_3$ group and there is a
surjective homomorphism from $H$ to $\mathbb Z$ with finitely generated
kernel. Then the FICwF is true for $\cal G$.\end{thm}

It is known that any nontrivial graph manifold has a finite sheeted cover
which either fibers over the circle or has large first Betti number (see
\cite{Lu}). We have already considered the case when $M$ is a virtual
fiber bundle. Hence we are left with the second case when there is a
finite sheeted cover with a large first Betti number.

Next we use the technique used in the proof of [\cite{R2}, main lemma] to
prove the FIC for another class of $3$-manifold groups.

Let $M$ be a closed $3$- manifold.

\begin{defn} \label{condition*}{\rm We say $M$ satisfies {\it Condition*}
if it
satisfies the following two conditions. 

\begin{itemize}

\item rank of $H_1(M, {\mathbb Z})$ is $\geq 1$,

\item either 

$(a).$ $[\pi_1(M), \pi_1(M)]$ is finitely generated or 

$(b).$ $M$ is a nontrivial graph manifold, rank of $H_1(M, {\mathbb Z})$
is greater than $1$, and $[\pi_1(M), \pi_1(M)]\cap \pi_1(P)$ is not a free
group for any Seifert fiber piece $P$ in the JSJT decomposition of $M$.
\end{itemize}}
\end{defn}

\begin{rem} \label{defrem} {\rm In the category of graph of groups, $(b)$
in Condition* has the following rewording. Note that $\pi_1(M)$ has the
structure of a graph of groups where each vertex represents a Seifert
fibered piece in $M$ and the associated group is the fundamental group of
the Seifert fibered piece. An edge is represented by a gluing
incompressible torus. As a subgroup of $\pi_1(M)$, $[\pi_1(M), \pi_1(M)]$
inherits a graph of group structure. $(b)$ in Condition* says that the
groups associated to each of the vertices of this graph are not free. Also
$(b)$ together with [\cite{He}, theorem 11.1] implies that $[\pi_1(M),
\pi_1(M)]$ is infinitely generated. This follows from the fact
that if the commutator subgroup is finitely generated and infinite cyclic
then $M$ is a Seifert fibered space which is a contradiction because $M$
is a nontrivial graph manifold. On the other hand if the commutator
subgroup is finitely generated and not infinite cyclic then $(2)$ of
[\cite{He}, theorem 11.1] implies that rank of $H_1(M, {\mathbb Z})$
is $1$ which is again a contradiction.

A class of examples which satisfy Condition* are described in
Section \ref{examples}}.\end{rem}

\begin{thm} \label{graph} Let $M$ be a closed $3$-manifold
such that there is a finite sheeted regular cover $\tilde M$ of $M$
satisfying Condition*. When $[\pi_1(\tilde M),\pi_1(\tilde M)]$ is
finitely generated assume that the FIC is true for $B$-groups. Then the 
FICwF is true for $\pi_1(M)$.\end{thm}

\begin{rem}{\rm Notice that when  $[\pi_1(\tilde M),\pi_1(\tilde M)]$ is
infinitely generated then we do not need the assumption that `the FIC is
true for $B$-groups'. Also when $[\pi_1(\tilde M),\pi_1(\tilde M)]$ is
finitely generated then one can deduce (using Proposition \ref{hempel})
that $\tilde M$ is finitely covered by a $3$-manifold fibering over the
circle.}\end{rem}

In this connection we pose the following problem.

\medskip
\noindent
{\bf Problem.} Let $M$ be a closed graph manifold. Then show that $M$
satisfies Condition*. 

Now let us recall that the fourth item of Theorem \ref{corollary} dealt
with $3$-manifolds
which are virtually of the type as in Stallings' fibration theorem. The
following result consider a similar but more general class of
$3$-manifold.

\begin{thm} \label{stallings} Let $M$ be a compact $3$-manifold so that
there is a finite index normal subgroup $G$ of $\pi_1(M)$ which sits in
the following exact sequence $$1\to H\to G\to Q\to 1$$ where $H$ is
finitely generated and $Q$ is an infinite group. If $H$ is not
infinite cyclic then assume that the FIC is true for any $B$-group. Then
the FICwF is true for $\pi_1(M)$.\end{thm} 

The following corollary is an immediate consequence of the results we
have proved.

\begin{cor} \label{coro} Let $G$ be any of the groups for which we have 
proved the FICwF, then the Whitehead group
$Wh(H)$, the reduced projective class group $\tilde K_0(H)$ and
$K_{-i}(H)$ for $i\geq 1$ vanish for any torsion free subgroup $H$ of
$G$.\end{cor}

The above corollary gives a partial solution to problem 3.32 of the R.
Kirby's problem list in \cite{K}. See [\cite{K}, p. 168] for a general
discussion about known results related to [\cite{K}, problem 3.32].
However, modulo Geometrization Conjecture, Corollary \ref{coro} is true
for any $3$-manifold $M$, whenever $\pi_1(M)$ is torsion free. Also recall
that $Wh(G)=\tilde K_0(G)=K_{-i}(G)=0$ for $i\leq 1$ is a classical result
of Waldhausen for a Haken $3$-manifold group $G$ (see \cite{W}) as
mentioned in the introduction.

\section{Proofs of the reduction theorem and its'
consequences}\label{reductionproof}

For notations and conventions used in this section we refer the reader to
\cite{DD}. The technique used are from the theory of action of groups on
trees and the main result we use is the structure theorem of groups acting
on trees.

\begin{proof}[Proof of the Reduction Theorem] Note that there is a
surjective homomorphism $\pi:G_1*G_2\to G_1\t G_2$ and by the
following Lemma the FIC is true for $G_1\t G_2$.

\begin{lemma} \label{lemma3} Let $G_1$ and $G_2$ be two groups and assume
that the FIC is true for both $G_1$ and $G_2$ then the FIC is true for 
the product $G_1\t G_2$.\end{lemma}

\begin{proof}  Consider the projection $p_1:G_1\t G_2\to
G_1$. By Lemma \ref{lemma C} we need to check that the FIC is true for
$p_1^{-1}(C)$ for any virtually cyclic subgroup $C$ of $G_1$. Note that
$p_1^{-1}(C)=C\t G_2$. Now consider the projection $p_2: C\t
G_2\to G_2$. Again we apply Lemma \ref{lemma C}. That is we need to show
that the FIC is true for $p_2^{-1}(C')$ for any virtually cyclic subgroup
$C'$
of $G_2$. But $p_2^{-1}(C')=C\t C'$ which is virtually poly-$\mathbb
Z$ and hence the FIC is true for $p_2^{-1}(C')$. This completes the proof
of
the Lemma.\end{proof}

Thus (by Lemma \ref{lemma C}) it is enough to check that the FIC is true
for $\pi^{-1}(C)$ for any
virtually cyclic subgroup $C$ of $G_1\t G_2$. We need the following Lemma.

\begin{lemma} \label{inverse} For a subgroup $C$ of $G_1\t G_2$ the
followings are true.
\begin{itemize}

\item $\pi^{-1}(C)$ is a countable free group if $C$ is either the trivial
group or the infinite cyclic group.

\item $\pi^{-1}(C)$ contains a countable free subgroup of finite index if 
$C$ is a finite group.

\end{itemize}

\end{lemma} 

\begin{proof} Let us first recall that the free product $G_1*G_2$ acts on
a tree $\cal T$ with vertex stabilizers $gG_1g^{-1}$ or $gG_2g^{-1}$, 
where $g$
varies in $G_1*G_2$, and trivial edge stabilizers. Hence any subgroup
$H$ of $G_1*G_2$ also acts on $\cal T$ with vertex stabilizers
$gG_1g^{-1}\cap H$ or $gG_2g^{-1}\cap H$, where $g\in G_1*G_2$. Now it is
easy to check
that $ker(\pi)\cap gG_ig^{-1}$ is trivial for $i=1,2$ and for any $g\in
G_1*G_2$. Hence using the structure theorem of groups acting on 
trees [\cite{DD}, I.3.4] it follows  that $ker(\pi)$ is a free group. Now
let $C$ be an infinite cyclic subgroup of $G_1\t G_2$. Then one 
checks that $\pi^{-1}(C)\cap gG_ig^{-1}$ is either infinite cyclic or the
trivial group for $i=1,2$ and for $g\in G_1*G_2$. Hence again we appeal to
[\cite{DD}, I.3.4] to deduce that $\pi^{-1}(C)$ is a free
group.

The second assertion is an immediate consequence of the first.

\end{proof} 

Let $C$ be a virtually cyclic subgroup of $G_1\t G_2$. If $C$ is
finite then $\pi^{-1}(C)<ker(\pi)\wr C$ and if $C$ is infinite and $C'$ is
an infinite cyclic normal subgroup of $C$ of finite index then
$\pi^{-1}(C)<\pi^{-1}(C')\wr (C/C')$. In both the cases, by Lemma
\ref{inverse} we get that $\pi^{-1}(C)$ is a subgroup of $F\wr H$ where
$F$ is a free group and $H$ is a finite group. The proof now follows
from the following lemma.

\begin{lemma}\label{fic} The FICwF is true for $F$ where $F$ is a
countable free group.\end{lemma}

\begin{proof} If $F$ is finitely generated then by Theorem 
\ref{Farrell-Jones}
the FICwF is true for $F$. If $F$ is infinitely generated then let
$F\simeq
\lim_{i\to\infty}F_i$ where $F_i$ is a finitely generated free group. We
get $F\wr H<\lim_{i\to\infty}(F_i\wr H)$. Now we appeal to Theorem
\ref{Farrell-Linnell} and Theorem \ref{Farrell-Jones} to deduce that the 
FICwF is
true
for $F$.\end{proof}

Using Lemma \ref{lemma A} we complete the proof of the Reduction
Theorem.\end{proof}   

\begin{proof}[Proof of Theorem \ref{refree}] Let $L$ be a finite
group. Note that there is a
surjective homomorphism $p:(G_1*G_2)\wr L\to (G_1\t G_2)\wr L$. It is
easy to see that $(G_1\t G_2)\wr L$ is a subgroup of $G_1\wr L\t
G_2\wr L$. By hypothesis and by Lemma \ref{lemma3} the FIC is true for
$(G_1\wr L)\t (G_2\wr L)$ and hence for $(G_1\t G_2)\wr L$ also.
Kernel of $p$ is $(\mbox{ker}(\pi))^L$, where $\pi$ is the surjective
homomorphism $G_1*G_2\to G_1\t G_2$, that is ker($p$) is a direct
product of $|L|$ copies of ker($\pi$). Let $C$ be a virtually cyclic 
subgroup of $(G_1\t G_2)\wr L$ and $C'$ be a cyclic subgroup
of $(G_1\t G_2)^L\cap C$ which is normal and of finite index in $C$. 
Let $C'$ be generated
by $\g=(\g_1,\cdots , \g_k)$ where $|L|=k$. Then the action of $\g$ on
$(G_1*G_2)^L$ is factorwise, that is the $i$-the coordinate of $\g$ acts
on the $i$-th factor of $(G_1*G_2)^L$. Without loss of generality we can 
assume that $\g_i\neq 1$ for $i=1,\ldots , l$ and $\g_i=1$ for
$i=l+1,\ldots , k$. Hence $p^{-1}(C')$ is a subgroup
of $\pi^{-1}(\l \g_1\r)\t\cdots\t\pi^{-1}(\l
\g_l\r)\t(\pi^{-1}(1))^{k-l}=H$(say). By Lemma 
\ref{inverse} each of the factor in the above expression is either 
free or contains a free subgroup of finite index. Now we have
$$p^{-1}(C)<p^{-1}(C')\wr
(C/C')<H\wr (C/C')$$$$<\pi^{-1}(\l \g_1\r)\wr (C/C') \t\cdots\t\pi^{-1}(\l
\g_l\r)\wr (C/C')\t(\pi^{-1}(1)\wr (C/C'))^{k-l}.$$ Using Lemma \ref{fic}
and the following easily verified Lemma 
we complete the proof of the Theorem.

\begin{lemma} \label{wr}Let $A$ and $B$ be finite groups and $G$ is any
group, then $(G\wr A)\wr B$ is a subgroup of $G^{A\t B}\wr (A\wr
B)$\end{lemma}\end{proof}

\begin{proof}[Proof of Corollary \ref{reductioncor}] At first note that 
the FICwF is true for finite groups and a compact irreducible
$3$-manifold is either aspherical or has finite fundamental group.
Therefore, for compact
$3$-manifolds the Corollary  follows from Lemma \ref{prime}, Theorem 
\ref{refree} and Lemma \ref{fic}. If the $3$-manifold is 
noncompact then we can write the manifold as an increasing union (under
inclusion) of compact submanifolds. Now the proof follows from the
previous case and by Theorem \ref{Farrell-Linnell}.\end{proof}

\begin{proof}[Proof of Corollary \ref{graphfic}] Let us first prove the
Corollary for finite graph of groups. Recall the notation $N_{\cal G}$
denoting the sum of number of vertices and number of edges of a finite
graph 
$\cal G$. The proof is by induction on $N_{\cal G}$. If $N_{\cal G}=1$
then there is only one vertex and no edge and hence FICwF is true for
$\pi_1(\cal G)$ by hypothesis. Assume that the Corollary is true
for finite graphs with $N_{\cal G}\leq n-1$. Let $\cal G$ be a graph with
$N_{\cal G}=n$. There are now two cases to consider.

\noindent
{\bf A.} Let $v$ be a vertex of $\cal G$ with only one edge
emanating from it. Let $G_v$ be the groups associated to the vertex $v$.
Then, clearly $\pi_1({\cal G})\simeq \pi_1({\cal G}_1)*G_v$ where ${\cal
G}_1$ is a graph with $N_{{\cal G}_1}=n-2$. Hence by the induction
hypothesis and by Theorem \ref{refree} FICwF is true for $\pi_1({\cal
G})$. 

\noindent
{\bf B.} All the vertices have valency $\geq 2$. Remove an edge from the
graph and let ${\cal G}_1$ be the resulting graph. Then
clearly $\pi_1({\cal G})\simeq \pi_1({\cal G}_1)*{\Bbb Z}$. Again by
induction hypothesis and Theorem \ref{refree} the proof is completed in
the finite graph case.

Now we consider the infinite graph case. Note that the fundamental group 
of any infinite graph which has countable number of vertices and edges can
be written as a direct limit of fundamental
groups of an increasing sequence of finite subgraphs. Using Theorem
\ref{Farrell-Linnell} and the previous case we complete the proof of the
Corollary.\end{proof} 

\section{Proofs of Theorems \ref{annulus}, \ref{geodesic} and the Main
Theorem}\label{pfmaintheorem} 

\begin{proof}[Proof of Theorem \ref{annulus}] Let $N$ be the double of
$M$. Then $N$ is a closed Haken $3$-manifold. If $N$ is a Seifert
fibered space then the FICwF is true for $\pi_1(N)$ by Theorem  
\ref{seifert}. Also since $\pi_1(M)$ is a subgroup of 
$\pi_1(N)$, FICwF is true for $\pi_1(M)$. So
assume that $N$ is not Seifert fibered. By Corollary \ref{jsjt} we 
can also assume that $N$ is not hyperbolic. Consider the
JSJT decomposition of $M$. This says that there is a collection of
finitely many
embedded incompressible tori (say $\cal T$), unique up to ambient isotopy,
in $M$ so that the complementary pieces are either Seifert fibered or
supports a complete hyperbolic metric. 

The proof of the Theorem will be
completed using the following two crucial lemmas.

\begin{lemma} \label{isotope} The family of tori ${\cal T}$ can be
isotoped in $N$ to be disjoint from $\p M\subset N$.\end{lemma}

\begin{proof} Assume on the contrary, that is, there is a member 
$T\in {\cal T}$ which can not be isotoped off $\p M$. 
Choose an isotope $T'$ of $T$ which intersects $\p M$ transversally so
that the number of circles in 
$C=\p M\cap T'$ is minimal. We can
assume that no circle in $C$ bounds a $2$-disc on $T'$. To see this note
that if there is such a disc then we can isotope $T'$ further to push it
off from $\p M$ making the number of circles in $C$ one less. Also
since $T$ cannot be isotoped off $\p M$ and since $N-\p M$ is disconnected
there are more than one circles
in $C$. Let $A$ be a component of $\overline {T'-C}$. Then $A$ is a
properly
embedded essential annulus in $M$. Which is a contradiction. This proves
the lemma.\end{proof}

\begin{lemma} \label{iso} $N$ is not a graph manifold.\end{lemma}

\begin{proof} At first by Lemma \ref{isotope} isotope $\cal T$ to make
it disjoint from $\p M$. Let $F$ be a component of $\p M$. Assume $N$ is a
graph manifold and choose a Seifert fibered piece $S$ in $N$ in the JSJT
decomposition of $N$ so that $F\subset S$. Since $N$ is not a Seifert
fibered space we conclude that $\p S\neq\emptyset$. Thus $S$ is a Seifert
fibered space with nonempty boundary and hence it has a finite sheeted
cover which is a product of a surface $F'$ and the circle. Note that $\p
F'\neq\emptyset$. Since $F$ is an incompressible surface in $S$ we get
$\pi_1(F)\cap (\pi_1(F')\t {\Bbb Z})$ is a finite index subgroup of
$\pi_1(F)$. On the other hand either $\pi_1(F)\cap (\pi_1(F')\t {\Bbb Z})$
is free or $\pi_1(F)$ contains a free abelian subgroup of rank $2$.
Which is a contradiction, since $F$ is a closed surface of genus $\geq 2$.
Hence $N$ is not a graph manifold.\end{proof}

The proof of the theorem is now easy. Since $N$ is a Haken $3$-manifold
and by Lemma \ref{iso} it is not a graph manifold, it follows that there
is a hyperbolic piece in the JSJT decomposition of $N$ and hence by
Corollary \ref{jsjt} FICwF is true for $\pi_1(N)$. Consequently for
$\pi_1(M)$ as well.

This completes the proof of the Theorem.\end{proof}

\begin{proof}[Proof of Theorem \ref{geodesic}] Let $M_1$ be a piece of
$M$ in the JSJT decomposition of $M$ so that $\p M\cap \p
M_1\neq \emptyset$. Then by Proposition \ref{agroup} the
interior of $M_1$ supports a complete hyperbolic metric. Since by 
hypothesis all the boundary components of $M$ are totally
geodesic, if $M_2$ is the double of $M_1$ along the boundary components of
genus $\geq 2$, then the interior of $M_2$ also supports a complete
hyperbolic Riemannian metric. Now let $N$ be the double of $M$. Then
$M_2\subset N$ and in the JSJT decomposition of $N$, $M_2$ is a hyperbolic
piece. Hence by \cite{Le}, $N$ supports a nonpositively curved Riemannian
metric. Consequently the FICwF is true for $\pi_1(N)$. By Lemma
\ref{lemma A} we complete the proof of the Theorem. 
\end{proof}

From the proof of Theorems \ref{annulus} and \ref{geodesic} we deduce the
following corollary.

\begin{cor} Let $M\in {\cal C}$ and $M$ does not contain any properly
embedded essential annulus or assume that each boundary component is
totally geodesic with respect to the metric as described in Proposition 
\ref{agroup}. Then $\pi_1(M)$ is a subgroup of the fundamental group of a
closed $3$-manifold P so that P is either Seifert fibered or is a
nonpositively curved Riemannian $3$-manifold.\end{cor}

\begin{proof}[Proof of the Main Theorem]
By hypothesis we have the following exact sequence. $$1\to K\to 
\pi_1(M)\to H\to 1.$$

\noindent
{\bf Case A.} Let us assume that $H$ is torsion free.

At first we check that the FIC is true for $\pi_1(M)$.

Let $M_L$ denote the covering of $M$ corresponding to a subgroup $L$ of
$\pi_1(M)$. When $L$ has infinite index then $M_L$ is a noncompact
$3$-manifold. Since the FIC is true for $H$, by
Lemma \ref{lemma C} we have
to check that the FIC is true for $p^{-1}(C)$ for any virtually cyclic
subgroup (including trivial) $C$ of $H$. Here $p$ is the homomorphism
$\pi_1(M)\to H$. Since $H$ is torsion free $C$ is either trivial or
infinite cyclic.

By the last assumptions in the statement of the theorem  
$M_{p^{-1}(C)}$ is a noncompact $3$-manifold. Hence by choosing a proper
smooth map
from $M_{p^{-1}(C)}$ to the real line we can write $M_{p^{-1}(C)}$ as a
union of increasing sequence of nonsimply connected, connected, compact
submanifolds $M_C^i$ with nonempty boundary. Hence $$p^{-1}(C)\simeq
\lim_{i\to \infty}\pi_1(M_C^i).$$ 

Since cover of an orientable irreducible
$3$-manifold is irreducible, $M_{p^{-1}(C)}$ is irreducible
(see theorem 3, section 7 of \cite{MY}). Therefore 
each $M_C^i$ has at least one boundary component of genus $\geq
1$. To see this we state and prove the following Lemma.

\begin{lemma}\label{boundary} Let $M$ be a nonsimply connected
irreducible $3$-manifold which is either noncompact or compact with
nonempty boundary. Let $N$ be a compact connected nonsimply connected
$3$-dimensional submanifold of $M$. Then there is at least one boundary
component of $N$ of
genus $\geq 1$.\end{lemma}

\begin{proof} On the contrary assume that all the boundary components of
$N$ are spheres, say $S_1,S_2,\cdots , S_k$. Since $M$ is irreducible there
is a $3$-disc $D_1\subset M$ so that $\p D_1=S_1$. Again since $M$ is
irreducible by [\cite{He}, lemma 3.8] $M-S_1$ is disconnected. Let
$M-S_1=M_1\cup int(D_1)$ where $M_1$ is connected. Since $N$ is connected
either $N\subset D_1$ or $N\subset \overline {M_1}$. In the former case it
is easy to show that $N$ is simply connected. Which is not possible by
hypothesis. So assume $N\subset
\overline {M_1}$. Let $N_1=N\cup D_1$. Then $N_1$ is a submanifold of $M$ 
with one less boundary components.
Continuing this process we find a closed submanifold (after capping off
all the boundary components of $N$) of $M$ of dimension $3$ - which is a
contradiction.\end{proof}  

By Theorem \ref{Farrell-Linnell} it is enough to check that the FIC
is true for $\pi_1(M_C^i)$. By capping off all the sphere boundary
components of $M_C^i$ we can assume that all the boundary components of 
$M_C^i$ are surfaces of genus $\geq 1$.

We need the following lemma.

\begin{lemma}\label{splitting} $\pi_1(M_C^i)\simeq G_1*\cdots *G_k*F^r$
where $F^r$ is a free group of rank $r$ and each $G_i$ is isomorphic to
the fundamental group of a compact irreducible $3$-manifold with
incompressible boundary.\end{lemma}

\begin{proof} By Lemma \ref{prime} $\pi_1(M_C^i)$ is isomorphic to a free
product of finitely many compact irreducible $3$-manifold (with
nonempty boundary of genus $\geq 1$ by Lemma \ref{boundary}) groups and a
free group of finite rank. We would like to use the Loop Theorem of
Papakyriakopolous to reduce the irreducible $3$-manifolds to irreducible
$3$-manifolds with nonempty incompressible boundary components.

We recall the Loop Theorem below.

\begin{thm3} ([\cite{He}, chapter 4, the loop theorem]) Let $M$ be a
compact $3$-manifold with boundary $\partial$. Let $\gamma$ be a
nontrivial element of $\pi_1(\partial)$ which goes to the trivial element
in $\pi_1(M)$. Then there is a simple closed curve on $\partial$ which
represents a nontrivial element in $\pi_1(\partial)$ and bounds a properly
embedded disc in $M$.\end{thm3} 

Let $N$ be a compact connected nonsimply connected irreducible
$3$-manifold with nonempty boundary. Then all the boundary components of
$N$ are surfaces of genus $\geq 1$.
Let $F$ be a boundary component of $N$ so that the kernel of $\pi_1(F)\to
\pi_1(N)$ is nontrivial. Then using the Loop theorem we get an properly
embedded disc $D$ in $N$. We cut $N$ along $D$ to get a new $3$-manifold
$N'\subset N$ with boundary. Using Lemma \ref{boundary} we get that
either both the components (if it is disconnected) of $N'$ are simply
connected or there is at least one boundary
component of genus $\geq 1$. If components of $N'$ are simply connected
then we stop
here. Otherwise we cap off all the sphere boundary components by
$3$-discs and denote it by $N'$ again. If $N'$ has two components
$N_1$ and $N_2$ then $\pi_1(N)\simeq \pi_1(N_1)*\pi_1(N_2)$. And if $N'$
has one component then $\pi_1(N)\simeq \pi_1(N')*{\Bbb Z}$. Now if some
boundary component of
$N'$ is compressible we apply again the Loop theorem. Standard technique
in $3$-manifold topology asserts that this process stops at finite stage.
That is after finitely many application of the Loop theorem we will get
finitely many $3$-manifolds which are either simply connected or
(after capping off sphere boundary components by $3$-discs) compact
$3$-manifolds with nonempty boundary and each boundary component is
incompressible. If any of these $3$-manifolds is not irreducible then we
apply Lemma \ref{prime} again.

This completes the proof of the Lemma.\end{proof}

By the Reduction Theorem and Theorem \ref{Farrell-Jones} we only have to 
check
that the FIC is true for $G_i$ for each $i$. So let $N$ be a compact
irreducible $3$-manifold with incompressible boundary components. If at
least one of the boundary component is a torus then by Corollary
\ref{jsjt} the FICwF is true for $\pi_1(N)$. So assume that all the
components of $\partial N$ are surfaces of
genus greater or equal to $2$. In this case by hypothesis the FICwF is
true for $\pi_1(N)$. Now the proof of the Theorem follows from Lemma 
\ref{splitting} and Theorem \ref{refree}  

This completes the proof of the Main Theorem in {\bf Case A}.

\noindent
{\bf Case B.} Let $J$ be a torsion free normal subgroup of $H$ of finite
index.
Then $\pi_1(M)$ is a subgroup of $p^{-1}(J)\wr (H/J)$. Note that
$p^{-1}(J)$ is again the fundamental group of a closed $3$-manifold so
that there is an exact sequence $1\to K\to p^{-1}(J)\to J\to 1$ satisfying
all the three properties in the statement of the Main Theorem and of
{\bf Case A}. 

Using Lemma \ref{wr} we see that to complete
the proof of the Main Theorem we only
have to check that the FICwF is true for $\pi_1(M^k)$ where $M$ is as in
the
statement with $H$ torsion free and $k$ is a
positive integer. 

Now we are in a more general setting as described in the following
two Lemmas.

\begin{lemma} \label{prfic} Assume that the FICwF is true for $H$. Then
the FICwF is
true for $H^k$ for any positive integer $k$.\end{lemma}

\begin{proof} The proof follows from Lemma \ref{lemma3} and by noting  
that $(H\t H)\wr G$ is a subgroup of $(H\wr G)\t (H\wr G)$. Here $G$ is a
finite group.
\end{proof}

\begin{lemma}\label{general} Let $K$ be a normal subgroup of a group $\G$
with quotient group $Q$ and $A:\G \to Q$ is the quotient homomorphism.
Assume the following.

\begin{itemize}
\item $A^{-1}(Z)\simeq \lim_{i\to\infty} \G ^Z_i$ where $Z$ is either
trivial or an infinite cyclic subgroup of $Q$ and $\{\G ^Z_i\}$ is a
directed system of groups so that for each $i$ the FICwF is true for
$\G ^Z_i$.

\item the FICwF is true for $Q$.

\end{itemize}

Then the FICwF is true for $\G$.\end{lemma}

Before we give the proof of Lemma \ref{general} let us first complete the
proof of the Main theorem using the Lemma. We only have to check that the
hypothesis of the above Lemma is satisfied. We have already seen that 
$p^{-1}(C)\simeq \lim_{i\to\infty}\pi_1(M_C^i)$ and $\pi_1(M_C^i)\simeq
G_1*\cdots *G_k*F^r$ where $F^r$ is a free group and each $G_i$ is the
fundamental group of a compact irreducible $3$-manifold with
incompressible boundary. Set $K=K, p^{-1}(J)=\G$ and $Q=J$. Then using 
Theorem \ref{refree}, Theorem \ref{Farrell-Jones},
Corollary \ref{jsjt} and {\bf Case A} we see that the
hypothesis of the Lemma is satisfied.

Thus we have completed the proof of the Main Theorem.\end{proof}

\begin{proof}[Proof of Lemma \ref{general}] Let $L$ be a finite group. We
have the following exact sequence. $$1\to K^L\to \G \wr L\to Q\wr L\to
1.$$
Here $K$ is the kernel of $A$. 

We denote the map $\G \wr L\to Q\wr L$ also by $A$.

By hypothesis the FIC is true for $Q\wr L$. Since $K\simeq
\lim_{i\to\infty}
\G ^{(1)}_i$ and the FIC is true for each $\G ^{(1)}_i$, by Theorem
\ref{Farrell-Linnell} the FIC is true for $K$ and hence for the product
$K^L$
by Lemma \ref{lemma3}. 

Let $Z'$ be a virtually cyclic subgroup of $Q\wr L$. If $Z'$ is finite
then $A^{-1}(Z') < K^L\wr Z' <\lim_{i\to\infty}((\G ^{(1)}_i\wr Z')\t
(\G ^{(1)}_i\wr Z')\t\cdots\t (\G ^{(1)}_i\wr Z'))$. There are $|L|$
number of factors inside the bracket. Hence applying Theorem
\ref{Farrell-Linnell}, Lemma \ref{lemma3} and Lemma \ref{lemma A} we get
that the FIC is true for $A^{-1}(Z')$. If $Z'$ is infinite let $Z$ be the
intersection of $Z'$ with the torsion free part of $Q^L<Q\wr L\simeq
Q^L\rtimes L$. Then $A^{-1}(Z')< A^{-1}(Z)\wr Z'/Z$. On the other hand
$A^{-1}(Z)\simeq K^L\rtimes Z$ where the action of $Z$ on $K^L$ is
factorwise. That is the $j$-th coordinate of a generator of $Z<Q^L$ acts
on the $j$-th factor of $K^L$. Hence $$K^L\rtimes Z < (K\rtimes
Z_1)\t\cdots\t (K\rtimes Z_{|L|}).$$ Here $Z_j$ is generated by
the $j$-th coordinate of a generator of $Z$. The right hand side of the
above expression is a subgroup of 
$$(\lim_{i\to\infty}\G ^{Z_1}_i)\t\cdots\t 
(\lim_{i\to\infty}\G ^{Z_{|L|}}_i)$$$$< 
\lim_{i\to\infty}(\G ^{Z_1}_i\t\cdots\t \G ^{Z_{|L|}}_i).$$

Hence $$A^{-1}(Z')<A^{-1}(Z)\wr
Z'/Z<\lim_{i\to\infty}(\G ^{Z_1}_i\t\cdots\t \G ^{Z_{|L|}}_i)\wr
Z'/Z$$$$<\lim_{i\to\infty}((\G ^{Z_1}_i\t\cdots\t \G ^{Z_{|L|}}_i)\wr
Z'/Z)$$$$<\lim_{i\to\infty}((\G ^{Z_1}_i\wr
Z'/Z)\t\cdots\t(\G ^{Z_{|L|}}_i\wr Z'/Z)).$$

By hypothesis the FIC is true for each of the factors inside the limit and
hence for $A^{-1}(Z')$ also. 

This proves the lemma.
\end{proof}

\section{Proof of Theorem \ref{corollary}}\label{corollaries}
\begin{proof} We give the proof of the theorem item wise.

\begin{itemize}
\item
At first we consider the
noncompact case. The compact case will follow from it. So let $M$ be a
noncompact $3$-manifold. Then by choosing a proper smooth map from $M$ to
the real line we can write $M$ as a union of
increasing sequence (under inclusion) of compact submanifolds $M_i$ with
nonempty boundary. By Lemmas \ref{boundary} and \ref{splitting}
$\pi_1(M_i)\simeq G_1*\cdots *G_k*F^r$ where each $G_i$ is isomorphic to
the fundamental group of a compact irreducible $3$-manifold with
incompressible boundary. The proof now follows by combining Theorem
\ref{Farrell-Linnell}, Theorem \ref{Farrell-Jones} and Corollary
\ref{jsjt}.

\item By the previous item we can assume that the manifold $M$ is closed.
Let $N$ be 
a finite sheeted cover of $M$ with first Betti number $\geq 2$. Then
$\pi_1(M)$ is a subgroup of $\pi_1(N)\wr G$ for some finite group $G$. If
we put $H=H_1(N, {\Bbb Z})$ then we check that all the three properties in
the statement of the Main Theorem are satisfied. Being a finitely 
generated abelian group, obviously $H$ has a
torsion free subgroup of finite index. $H\wr G$ is a virtually poly-${\Bbb
Z}$ group and hence the FIC is true for $H\wr G$ (Proposition
\ref{prop2.4}) for any finite group $G$.
Since $H$ has rank greater or equal to $2$ the third condition is
satisfied. Hence by the Main Theorem the FICwF is true for
$\pi_1(N)$.

\item Let $M$ be a
nontrivial graph manifold. That is there is an incompressible torus
embedded in the manifold. Now we appeal to \cite{Lu} where it is proved
that in this situation either $M$ has a finite sheeted cover which is a
torus bundle over the circle or there is a finite sheeted cover of $M$
with arbitrarily large first Betti number. In the first case it follows
that $\pi_1(M)$ is virtually poly-${\Bbb Z}$ for which the FIC is true
(Proposition \ref{prop2.4}). In the second case the result follows from
the last item. 

For the second assertion recall that Thurston's Geometrization
Conjecture claims that any compact irreducible $3$-manifold 
either has a finite fundamental group or is Seifert fibered, Haken or
hyperbolic. Using Theorem \ref{seifert} we complete the 
proof of this item.

\item By [\cite{FR}, algebraic lemma] and Lemma \ref{prfic} it is enough
to prove
the FICwF for any fibered $3$-manifold. So let $M$ be a
$3$-manifold fibering over the circle. Let $F$ be the fiber of the fiber
bundle projection $M\to {\Bbb S}^1$. If $F$ is the sphere, torus or the
Klein bottle then $\pi_1(M)\wr G$ is virtually poly-${\Bbb Z}$ and hence
the proof is complete in this case. So assume the $F$ has genus $\geq 2$. 

Since the commutator subgroup is characteristic we have the following
exact sequence. $$1\to [\pi_1(F), \pi_1(F)]\to \pi_1(M)\to H_1(F, {\Bbb
Z})\rtimes \langle t\rangle \to 1.$$ 

Let $H=H_1(F, {\Bbb Z})\rtimes \langle t\rangle$. Then it is easy to
check all the three properties in the statement of the Main Theorem for
the group $H$. 

\item The proof of this case follows from the proof of the third item.
\end{itemize}
\end{proof}

\section{Proofs of Theorems \ref{seifert}, \ref{pdgroup}, \ref{graph} and
\ref{stallings}}\label{manifoldproof}

\begin{proof} [Proof of Theorem \ref{seifert}] If $\pi_1(S)$ is finite
then it is known that
the FICwF is true for this group. So assume that $\pi_1(S)$ is infinite.
Then
there exists the following well-known exact sequence. $$1\to C\to
\pi_1(S)\to \pi_1^{orb}(B)\to 1.$$ Here $\pi_1^{orb}(B)$ denotes the
orbifold fundamental group of the base orbifold of the Seifert fibered
space $S$ and $C$ is an infinite cyclic group generated by a regular
fiber of the Seifert fibration of $S$. Let $L$ be a finite group and
consider the following exact sequence. $$1\to C^L\to
\pi_1(S)\wr L\to \pi_1^{orb}(B)\wr L\to 1.$$ Note that when $B$ is
closed, $\pi_1^{orb}(B)$ contains a normal closed surface subgroup of
finite index and hence by Theorem \ref{Farrell-Jones} and Lemma 
\ref{prfic} the
FICwF is
true for $\pi_1^{orb}(B)$. In the compact with nonempty
boundary case take the double of $B$ and apply Lemma \ref{lemma A}. Let
$p:\pi_1(S)\wr L\to \pi_1^{orb}(B)\wr L$ be the projection map and $V$ is
a virtually cyclic subgroup of $\pi_1^{orb}(B)\wr L$. Then $p^{-1}(V)$ is 
virtually poly-${\mathbb Z}$ and hence the FIC is true for $p^{-1}(V)$ by
Proposition \ref{prop2.4}. Thus by Lemma \ref{lemma C} the FIC is true
for $\pi_1(S)\wr L$.\end{proof}

\begin{proof} [Proof of Theorem \ref{pdgroup}] If the kernel $K$ of the
homomorphism $H\to
\mathbb Z$ is finite then $\cal G$ is virtually cyclic and hence the FIC
is true
for $\cal G$. So we assume that $K$ is finitely generated and infinite. By
[\cite{Hi}, theorem 1.20] $K$ is a $PD_2$ group and hence by [\cite{EL}, 
\cite{EM}] $K$ is isomorphic to the fundamental group of a closed
surface of positive genus. Since any automorphism of a closed surface is
induced by a diffeomorphism of the surface it follows that $H$ is
isomorphic to the fundamental group of a closed $3$-manifold $L$ which
fibers over ${\mathbb S}^1$. Thus ${\cal G}<\pi_1(L)\wr U$ where $U$ is a
finite group. By Theorem \ref{corollary} we complete the
argument.\end{proof}

We need the following lemma to prove Theorem \ref{graph}. 

\begin{lemma} \label{cruciallemma}Let $M$ be a compact graph
manifold and $\Gamma$ is a normal subgroup of $\pi_1(M)$ of infinite index
and with the property that $\Gamma\cap \pi_1(P)$ is not free for any
Seifert fibered piece $P$ in the JSJT decomposition of $M$. Let
$M_\Gamma$ be the covering of $M$ corresponding to $\Gamma$. Then
$M_\Gamma=\cup_iN_i$ where for each $i=1,2,\cdots,$ $N_i\subset N_{i+1}$
and $N_i$ is a compact irreducible $3$-manifold with incompressible tori
boundary components.\end{lemma}

\begin{proof} Let $f:M_\Gamma\to M$ be
the covering projection and let $P$
be a Seifert fibered piece of $M$. Let $\tilde P$ be a component of
$f^{-1}(P)$. We will show that $\tilde P$ is a Seifert fibered space.  
Note that $M_\Gamma$ is noncompact and has empty boundary. 
Noncompactness implies there is a smooth proper surjective function
$\delta: M_\Gamma\to {\mathbb R}^{\geq 0}$ or $\delta: M_\Gamma\to
{\mathbb R}$ according as $M_\Gamma$ has $1$ or more than $1$ ends
respectively.  

The proof of the lemma now consists of similar ideas as
in the proof of [\cite{R2}, corollary 6.2]. At first recall that there is
a central infinite cyclic subgroup $Z$ of $\pi_1(P)$ with quotient 
$\pi_1^{orb}(B_P)$, where $B_P$ is the base orbifold of $P$ (\cite{He}).
We claim that $Z\cap (\Gamma\cap \pi_1(P))$ is nontrivial. If not then
$\Gamma\cap\pi_1(P)$ injects into $\pi_1^{orb}(B_P)$. Since $B_P$ has
nonempty boundary and $\Gamma\cap\pi_1(P)$ is torsion free we get that
$\Gamma\cap\pi_1(P)$ is free - which is a contradiction. Now since
$f_*(\pi_1(\tilde P))$ is conjugate to $\Gamma\cap \pi_1(P)$,
$\pi_1(\tilde P)$ contains an infinite cyclic central subgroup. Hence we
can apply [\cite{Ma}, theorem 1.1] to deduce that $\tilde P$ is Seifert
fibered. 

Now we proceed to give a filtration of $\tilde P$ in an
increasing union of compact irreducible $3$-manifolds with
incompressible tori boundary components. Let $B_{\tilde P}$ be the base
surface of $\tilde P$. Choose
a filtration of $B_{\tilde P}$ by increasing union of compact subsurfaces 
with incompressible circular boundary components so that the boundary
components do not contain any of the orbifold points. Now pulling back
this filtration by the quotient map $\tilde P\to B_{\tilde P}$ we get a
filtration of $\tilde P$ by compact irreducible $3$-manifolds $\tilde P_i$
with incompressible tori boundary components. Thus we get a covering of
$M_\Gamma$ by
the collection $\{\tilde P_i\}_{i,P}$ where $i\in {\mathbb N}$ and $P$
varies over Seifert fibered pieces of $M$ and $\tilde P$ varies over all
components of $f^{-1}(P)$. Notice that any two members of the covering
$\{\tilde P_i\}_{i,P}$ either do not intersect or intersect along some
tori boundary or along some incompressible annuli on tori boundary
components.
Hence union of any finitely many members of this covering is a compact
irreducible (possibly not connected) $3$-manifold with incompressible
tori boundary components. 

Now we can choose $0< r_1<r_2\cdots <r_m<\cdots$ and $i^l_1,i^l_2,\cdots ,
i^l_{k_l}$, $l=1,2,\cdots$, so that $r_m\to\infty$ and
$$\cup_{j=1}^{j=k_{l-1}} {\tilde P}_{i^{l-1}_j}\subset \delta^{-1}([-r_l,
r_l])\subset \cup_{j=1}^{j=k_l} {\tilde P}_{i^l_j}$$ and
$\cup_{j=1}^{j=k_l} {\tilde P}_{i^l_j}$ is connected. Write
$N_l=\cup_{j=1}^{j=k_l} {\tilde P}_{i^l_j}$.

This completes the proof of the Lemma.
\end{proof}

\begin{proof} [Proof of Theorem \ref{graph}] Let $M$ be a closed 
$3$-manifold which satisfies Condition* virtually. That is there is a
regular finite sheeted covering $\tilde M$ of $M$
which satisfies Condition*. Let $L$ be the group of covering
transformation of $\tilde M\to M$. We have $\pi_1(M)<\pi_1(\tilde M)\wr
L$. Thus we only need to check the FIC for $\pi_1(\tilde M)\wr L$. There
are two cases to consider. 

\medskip
\noindent
{\bf Case A.} $[\pi_1(\tilde M), \pi_1(\tilde M)]$ is finitely generated.
Since $H_1(\tilde M, {\mathbb Z})$ is infinite we can apply Proposition
\ref{hempel} to deduce that $[\pi_1(\tilde M), \pi_1(\tilde M)]$  is
isomorphic to the fundamental group of a compact surface (say $F$). If
$[\pi_1(\tilde M), \pi_1(\tilde M)]$ is infinite cyclic then $\pi_1(M)$
is virtually poly-${\mathbb Z}$ and hence the FICwF is true for $\pi_1(M)$
by Proposition \ref{prop2.4}.
Otherwise we apply Proposition \ref{hempel} again to conclude that
$H_1(\tilde M, {\mathbb Z})$ is virtually infinite cyclic. Let $C$ be an
infinite cyclic normal subgroup of $H_1(\tilde M, {\mathbb Z})$ of finite
index and let $A:\pi_1(\tilde M)\to H_1(\tilde M, {\mathbb Z})$ denotes 
the
abelianization homomorphism. Hence $A^{-1}(C)$ is a
finite index subgroup of $\pi_1(\tilde M)$
and thus isomorphic to the fundamental group of a closed $3$-manifold. 
Also since $A^{-1}(C)$ has a finitely generated normal subgroup with 
infinite cyclic quotient, by Stallings' fibration theorem 
$A^{-1}(C)$ is isomorphic to the fundamental group of a
closed $3$-manifolds fibering over the circle. 

Let $C'=H_1(\tilde M, {\mathbb Z})$. Then $A^{-1}(C')< A^{-1}(C)\wr
(C'/C)$. Hence by Lemma \ref{wr} $\pi_1(\tilde M)\wr
L=( A^{-1}(C'))\wr L <(A^{-1}(C)\wr (C'/C))\wr L < (
A^{-1}(C))^{((C'/C)\t L)}\wr ((C'/C)\wr L)$. By Theorem \ref{corollary} 
(item $4$) and Lemma \ref{lemma3} the FIC is true for 
$(A^{-1}(C))^{((C'/C)\t L)}\wr
((C'/C)\wr L)$ and hence for $\pi_1(\tilde M)\wr L$ also. This completes
the proof of the theorem in this case.

\noindent
{\bf Case B.} $M$ satisfies $(b)$ of Condition*. By Remark \ref{defrem}
$[\pi_1(\tilde M), \pi_1(\tilde M)]$ is infinitely generated. 

By Lemma \ref{cruciallemma} it follows that $A^{-1}(Z)< \lim_i G^Z_i$
where each
$G^Z_i$ is isomorphic to the fundamental group of a compact irreducible 
$3$-manifold with incompressible tori boundary components and $Z$ is any
infinite (or trivial) cyclic
subgroup of $H_1(\tilde M, {\mathbb Z})$ and $A$ denotes the homomorphism 
$\pi_1(\tilde M)\to H_1(\tilde M, {\mathbb Z})$. By Corollary \ref{jsjt} 
the FIC is true for $G^Z_i\wr L$ for any finite group $L$. 

Now note that we are in the setting of the Lemma \ref{general}.

This completes the proof of Theorem \ref{graph}.
\end{proof}

\begin{proof}[Proof of Theorem \ref{stallings}] Let us first check that
the hypothesis of the theorem implies that $N$ is irreducible where $N$
is a finite sheeted covering of $M$ realizing $G$. That is $G=\pi_1(N)$.
If $N$ is not irreducible, then either $G$ is a nontrivial free product 
or $N$ is an ${\mathbb S}^2$-bundle over the circle. Since $H$ is
nontrivial and of infinite index in $G$ the second case does not occur.
Since $H$ is a nontrivial finitely generated normal subgroup of the free
product $G$, by [\cite{He}, lemma 11.2] $H$ should have finite index in
$G$. Again this is a contradiction. Hence $N$ is irreducible. 

The main ingredient behind the proof is [\cite{He}, theorem 11.1]. We
need the following proposition which is an application of [\cite{He},
theorems 11.1, 11.6] and the proof of some special cases of the Seifert
fiber space conjecture. 

\begin{prop} \label{hempel} Let $N$ be a compact $3$-manifold so that
there is an exact sequence $1\to H\to \pi_1(N)\to Q\to 1$, where $H$ is
finitely generated and $Q$ is infinite. Then the followings hold.

\begin{itemize}

\item $H$ is the fundamental group of a compact surface.

\item if $H\simeq {\mathbb Z}$ then either $N$ is Seifert fibered or
there is a finite sheeted Seifert fibered covering of $N$.

\item if $H$ is not infinite cyclic then $Q$ is virtually cyclic.

\item if $Q$ is cyclic then either ${\cal P}(N)$ fibers over the circle
with fiber a compact surface with fundamental group isomorphic to $H$ or
${\cal P}(N)$ is homotopy equivalent to ${\mathbb P}^2\t {\mathbb
S}^1$.
\end{itemize}  
\end{prop}

\begin{rem} \label{poincare} {\rm For definition of the Poincar\'{e}
associate ${\cal P}(N)$ of $N$ see [\cite{He}, p. 88]. For our purpose it
is enough to recall that $N$ and ${\cal P}(N)$ have the same fundamental
group.}\end{rem}

\begin{proof} [Proof of Proposition \ref{hempel}] The first and third
conclusions are $(1)$ and $(2)$ of
[\cite{He}, theorem 11.1] respectively. For the second one note that
$N$ is a compact irreducible $3$-manifold with an infinite cyclic normal
subgroup $H$. If $H$ is central and $N$ is orientable then by a by now
well-known theorem $N$ is Seifert fibered. Otherwise by [\cite{He}, remark
11.2] there is a two sheeted covering $p:\tilde N\to N$ so that $H\cap
p_*(\pi_1(\tilde N))$ is
central in $p_*(\pi_1(\tilde N))$. Hence using the previous case we
conclude the proof of the second item. The fourth item is $(1)$ of 
[\cite{He}, theorem 11.6]. This completes the proof of the Proposition.
\end{proof}

To prove the Theorem there are now two cases to consider.

\noindent
{\bf Case A.} $H$ is infinite cyclic. Since $H$ is now an infinite cyclic
normal subgroup of the fundamental group of the compact irreducible
$3$-manifold $N$, by Proposition \ref{hempel} either $N$ is Seifert
fibered or there is a finite sheeted cover of $N$ which is Seifert
fibered. Hence there is a finite sheeted covering of $M$ which is 
Seifert fibered. Therefore the FICwF is true for $\pi_1(M)$ by Theorem 
\ref{seifert}. 

\noindent
{\bf Case B.} $H$ is not infinite cyclic. In this case by the third case
of Proposition \ref{hempel}, $Q$ is virtually cyclic. Let $Q'$ be an
infinite cyclic normal subgroup of $Q$ of finite index. Let $p:G\to Q$ be
the projection. Then applying the fourth case of Proposition
\ref{hempel} to the exact sequence $1\to H\to p^{-1}(Q')\to Q'\to 1$ we
get that ${\cal P}(M')$ either fibers over the circle or ${\cal P}(M')$ is
homotopically equivalent to ${\mathbb P}^2\t {\mathbb S}^1$. Here $M'$
is the finite sheeted covering of $N$ corresponding to the subgroup
$p^{-1}(Q')$. Hence we get that either ${\cal P}(M)$ has a finite sheeted
covering which fibers over the circle or $\pi_1(M)$ is virtually infinite
cyclic. Using Theorem \ref{corollary} we complete the proof of the
theorem.\end{proof}

\section{Properties of the class $\cal C$}\label{bgroup}

This section is devoted to describe some properties of $3$-manifolds
belonging to the class $\cal C$. Let ${\cal C}'$ be the subclass of  $\cal
C$ consisting of $3$-manifolds which do not contain any properly embedded
essential annulus. We establish that, in some suitable sense (defined
below), ${\cal C}'$ is a basis of $\cal C$, that is $\cal C$ can be
obtained from ${\cal C}'$ by some fundamental operations.

Let $M$ be a compact connected $3$-manifold with nonempty boundary. Let
$A_1$ and $A_2$ be two disjoint incompressible annuli on $\p M$. Let us
denote by $M\#_A$
the manifold obtained by identifying $A_1$ with $A_2$ by some
diffeomorphism of $A_1$ and $A_2$. Now let $M_1$ and $M_2$ be two compact
connected $3$-manifolds with nonempty boundary. Let $A_1\subset \p M_1$
and $A_2\subset \p M_2$ be two incompressible annuli. We denote by
$M_1\#_A M_2$
the manifold obtained by identifying $A_1$ and $A_2$ by some
diffeomorphism. On the fundamental group level the first case corresponds
to $HNN$-extension and the second one corresponds to generalized free
product case.

\begin{defn} {\rm We call $M\#_A$ an $HNN-${\it connected sum} along an  
annulus and $M_1\#_A M_2$ a {\it generalized connected sum} along an
annulus.}\end{defn}

\begin{prop}\label{building} Let $M\in {\cal C}$. Then there are  
$M_1,\ldots , M_k\in {\cal C}'$  
and $M$ is obtained from $M_1,\ldots , M_k$ by successive
application of $HNN$ and generalized connected sum along
annuli on the boundary components of $M_i$.\end{prop}

\begin{proof} If $M$ has no properly embedded essential annulus then there
is nothing to prove. So assume there is a properly embedded essential
annulus $A\subset M$. Let $N=\overline {M-A}$. Then $N$ is an irreducible
$3$-manifold and it is easy to show by Euler characteristic calculation
that all the boundary components of $N$ are surfaces of genus $\geq 2$.
Incompressibility of the boundary also follows easily. Hence $N\in {\cal
C}$. If $N$ is connected then $M$ is obtained from $N$ be $HNN-$ connected
sum along an annulus and if $N$ is disconnected then $M$ is obtained from
$N$ by generalized connected sum along an annulus. Now again we apply the
same procedure if $N$ contains a properly embedded essential annulus.
Applying standard (reducing complexity) tricks from $3$-manifold topology
it follows now that after a finite stage we get manifolds which contain
no properly embedded essential annuli. This proves the 
proposition.\end{proof} 

\begin{rem}{\rm Note that the splitting theorem in [\cite{JS}, p. 157]
states that any compact Haken $3$-manifold contains a two-sided
incompressible $2$-manifold which is unique up to ambient isotopy and the
components of this $2$-manifolds are either tori or annuli and none of
them is boundary parallel. Also after cutting the manifold along this
$2$-manifold one gets pieces which are either Seifert fibered or simple.
Here we would like to remark that the proof of the Proposition 
\ref{building} consists of only cutting along these unique class of
annuli.}\end{rem}

In the following we give examples of members of $\cal C$ which contains 
essential annuli and the boundary is not totally geodesic. 

\begin{exm} \label{examp} {\rm Let $S$ be a closed orientable surface of
genus $\geq 2$.
Let $M=S\t [0,1]$. Then $M\in {\cal C}$. It is also easy to show that the
interior of $M$ supports a complete hyperbolic metric. We check the
following two properties about $M$.

\begin{itemize}
\item $M$ contains properly embedded essential annulus.
\item at least one of the boundary components of $M$ is not totally
geodesic.
\end{itemize}

The first property is obvious. For example take any simple closed curve
$\gamma\subset S$ and then consider $\gamma\t [0,1]$. 

For the second property, suppose on the contrary, that the boundary
components are totally geodesic. Identify $S\t \{0\}$ with $S\t \{1\}$ by
the identity map and let $N$ be the resulting manifold. Since $S\t \{0\}$
and $S\t \{1\}$ are totally geodesic in $M$, $N$ also supports a
hyperbolic metric. This is a contradiction because by Preissman's theorem
([\cite{C}, theorem 3.2]) any nontrivial abelian subgroup of $\pi_1(N)$
is infinite cyclic, on the other hand $\pi_1(N)$ contains a free abelian
subgroup on $2$ generators.}\end{exm}

\section{On $A$-regular Riemannian manifolds and
$3$-manifolds}\label{arrm}

This section is devoted to give some examples of noncompact Riemannian
manifolds with $A$-regular Riemannian metric. The basic material on
Riemannian geometry we use here are given in any standard book (for
example see \cite{Hic} or \cite{C}). We follow the notations and
conventions of \cite{C}. 

\begin{prop}\label{aregular1} Let $M$ be a $A$-regular
Riemannian manifold. Then the product metric on the $n$-fold product $M^n$
is also $A$-regular.\end{prop}

\begin{proof} It is enough to prove the Proposition for $n=2$. Let
$\chi(K)$ be the space of all vector fields on a
manifold $K$ and let $R_K(X,Y,Z,W)=\l R_K(X,Y)Z, W\r$, $X,Y,Z,W\in\chi(K)$ 
denote the curvature tensor if $K$ is a Riemannian manifold. Here $\l\ , \ 
\r$ denotes the Riemannian metric on $K$ and $R_K(X,Y)Z$ is the curvature.
Let $p_1:M\t M\to M$ and
$p_2:M\t M\to M$ be the first and the second projection. For any vector
field $V$ on $M\t M$ let $V_1=(p_1)_*(V)$ and $V_2=(p_2)_*(V)$ where
$(p_i)_*:T(M\t M)\to TM$ is the induced map on the tangent bundle for each
$i=1,2$. Let $X,Y,Z,W\in\chi(M\t M)$ and $X=X_1+X_2, Y=Y_1+Y_2, Z=Z_1+Z_2,
W=W_1+W_2$ be the orthogonal decomposition of the vector fields
as defined above.
where $X_i,Y_i,Z_i,W_i\in\chi(M)$ for $i=1,2$. 

Now we need the following easily verified lemma.

\begin{lemma} \label{connection}Let $M$ and $N$ be two Riemannian
manifolds and $p_1:M\t N\to M$ and $p_2:M\t N\to N$ be the two
projections. Let $X,Y,Z$ and $W$ be vector fields on $M\t N$ such
that $(p_2)_*(X)=0=(p_1)_*(Y)$ and $(p_2)_*(W)=0=(p_1)_*(Z)$ and let
$\nabla$ denotes the Riemannian connection on $M\t N$. Then the followings
are true.

\begin{itemize}
\item $\nabla_XY=\nabla_WZ=[X,Y]=[W,Z]=0$.
\item let $f:N\to {\Bbb R}$ be a smooth function then $X(f\circ p_2)=0$.
\end{itemize}
\end{lemma}

Now recall the definition of the curvature. For a Riemannian manifold
$K$, $$R_K(X,Y)Z=\nabla_Y\nabla_XZ-\nabla_X\nabla_YZ+\nabla_{[X,Y]}Z.$$

Expanding the curvature linearly and using the above Lemma we deduce 
that $$R_{M\t M}(X,Y,Z,W)=R_{M\t M}(X_1,Y_1,Z_1,W_1)+R_{M\t
M}(X_2,Y_2,Z_2,W_2).$$ Fix $X,Y,Z$ and $W$. Let $(p,q)\in M\t M$. Then the
two real numbers $$R_{M\t M}(X_1(p,q),Y_1(p,q),Z_1(p,q),W_1(p,q))$$ and
$$R_M(X_1(p),Y_1(p),Z_1(p),W_1(p))$$ are equal and similarly $$R_{M\t
M}(X_2(p,q),Y_2(p,q),Z_2(p,q),W_2(p,q))$$ and 
$$R_M(X_2(q),Y_2(q),Z_2(q),W_2(q))$$ are equal. 

Hence
$$R_{M\t M}(X^1(p,q),\cdots ,
X^4(p,q))$$$$=R_M(X^1_1(p),\cdots ,
X^4_1(p))+R_M(X^1_2(q),\cdots ,X^4_2(q)).$$

Now using the above Lemma and by linearity one easily checks the 
following equality. $$\nabla^iR_{M\t
M}(X^1(p,q), X^2(p,q),\cdots ,
X^{i+4}(p,q))$$$$=\nabla^iR_M(X^1_1(p),\cdots ,
X^{i+4}_1(p))+\nabla^iR_M(X^1_2(q),\cdots ,X^{i+4}_2(q))$$ for each $i$.
Here for each $j=1,2,\cdots, i+4$, $X^j_1$ and $X^j_2$ denotes the two
orthogonal components of $X^j$,
that is $X^j=X^j_1+X^j_2$.

Thus if $M$ is an $A$-regular Riemannian manifold and $A_i$ are
nonnegative constants such that $|\nabla^iR_M|\leq A_i$ for each $i$ then
$|\nabla^iR_{M\t M}|\leq 2A_i$ for each $i$. Hence $M\t M$ is also
$A$-regular.

This proves the Proposition.\end{proof}

In the proof of the above proposition, in fact we have shown something
more.

\begin{cor} Let $M$ and $N$ be two complete Riemannian manifolds with
$A$-regular metric. Then the product metric on $M\t N$ is also
$A$-regular.\end{cor}

\begin{prop}\label{aregular} Let $M$ be a complete Riemannian manifold
and outside a compact subset it has constant sectional curvature. Then
the metric is $A$-regular.\end{prop}

\begin{proof} Note that the curvature tensor is completely determined by
the sectional curvatures ([\cite{C}, chapter 4, lemma 3.3]). It is now
easy to check that the covariant derivatives of the curvature tensor is
zero outside a compact subset. And obviously they are bounded on a compact
subset. This proves the Proposition.\end{proof} 

An application of Proposition \ref{aregular} is Proposition \ref{agroup}.

\begin{proof} [Proof of Proposition \ref{agroup}] Note that $M$ is a Haken 
$3$-manifold and hence it admits a
JSJT decomposition (along tori) into Seifert fibered and hyperbolic
pieces. Since
boundary components of Seifert fibered pieces are tori the pieces abutting
the boundary components of $M$ are all hyperbolic. This proves the first 
assertion. On the other hand
by \cite{Le} the hyperbolic metric in the hyperbolic pieces can be
deformed near the tori boundary components and a Riemannian metric of
nonpositive sectional curvature can be given in the interior of $M$ so
that outside a compact subset of $M$ the metric has constant $-1$
sectional curvature. Hence this Riemannian metric on $M$ is $A$-regular
by Proposition \ref{aregular}. Also if we consider a finite product
$M\t\cdots\t M$ then the product Riemannian metric is also $A$-regular and 
nonpositively curved by Proposition \ref{aregular1}. Now since
$\pi_1(M)\wr G\simeq \pi_1(M^G)\rtimes G$ it follows that $\pi_1(M)\wr G$
is an $A$-group.\end{proof} 

\section{Examples of $3$-manifolds satisfying Condition*}\label{examples} 

In this section we work out some examples of graph manifolds for which
Condition* is satisfied. The simplest example of this sort is given in the
following way.

\begin{exm} {\rm Let $F_1$ and $F_2$ be two compact orientable non-simply
connected surfaces
and $P_i=F_i\t
{\Bbb S}^1$ for $i=1,2$. Assume $F_i$ has one boundary
component for $i=1,2$. Let for each $i=1,2$, 
$(\lambda_i,\mu_i)$ denotes the basis of $\pi_1(\partial P_i)\simeq {\Bbb
Z}\t {\Bbb Z}$ where $\mu_i$
represents $\partial F_i$ and $\lambda_i$ represents the second factor
${\Bbb S}^1$. Since $\p F_i$ has one component, $\mu_i$ represents an
element of the commutator subgroup of $\pi_1(F_i)$
for $i=1,2$. 

Let $V$ be the closed manifold obtained by identifying $P_1$ and $P_2$
along their boundary tori by a diffeomorphism $f$. Then $V$ is a graph
manifold. 

The diffeomorphism $f$ induces the following isomorphism $$f_*:\pi_1(\p
P_1)\simeq {\Bbb Z}\t {\Bbb Z}\to \pi_1(\p P_2)\simeq {\Bbb Z}\t {\Bbb
Z}$$
and hence there are pairs of integers $(p_1, q_1)$ and $(p_2, q_2)$ with
the following properties.

\begin{itemize}
\item $p_1q_2-p_2q_1=1$
\item $f_*$ sends $\lambda_1$ to $\lambda_2^{p_1}\mu_2^{q_1}$ and $\mu_1$
to $\lambda_2^{p_2}\mu_2^{q_2}$. 
\end{itemize}

\begin{prop} Under the above notations the followings are true. 

\begin{itemize}
\item if $p_2\neq 0$ then $V$ satisfies $(b)$ in the definition of
Condition*.
\item if $p_2=0=q_1$ then $V$ is a product.
\end{itemize}

\end{prop}

\begin{proof} By an easy calculation we see that in the fundamental group
of $V$ the following are satisfied.

\begin{itemize}
\item $\lambda_1^{p_2}=\mu_1^{p_1}\mu_2^{-1}$.
\item $\lambda_2^{p_2}=\mu_1\mu_2^{-q_2}$.
\end{itemize}

If $q_1=0=p_2$ then it easily follows that $V$ is a product of a
compact surface and the circle. 

Now if $p_2\neq 0$ then we have that 
$\lambda_2^{p_2}=\mu_1\mu_2^{-q_2}$ and
$\lambda_1^{p_2}=\mu_1^{p_1}\mu_2^{-1}$ in $\pi_1(V)$. Since $\mu_1$ and 
$\mu_2$ both represents elements of $[\pi_1(V),\pi_1(V)]$ we see that
$[\pi_1(V),\pi_1(V)]\cap \pi_1(P_2)$ contains $\lambda_2^{p_2}$ and hence
$[\pi_1(V),\pi_1(V)]\cap \pi_1(P_2)$ is not free. This follows from the
fact that $[\pi_1(V),\pi_1(V)]\cap [\pi_1(F_2), \pi_1(F_2)]$ is
nontrivial. Similarly $[\pi_1(V),\pi_1(V)\cap \pi_1(P_1)$ is also not
free.

This proves the Proposition.\end{proof}

Similar ideas as in the above Proposition can be used to construct
examples satisfying $(b)$ of Condition*.}\end{exm}

\begin{rem}{\rm Here we remark that there are only
finitely many ways to identify the torus boundary of $P_1$ with
that of $P_2$ to produce a closed graph manifold admitting a metric
of nonpositive sectional curvature (see \cite{KL}).}\end{rem} 

\medskip
\noindent
{\bf Acknowledgment.} The author would like to thank F.T. Farrell for
pointing out the relevance of the arguments used in the proof of
[\cite{FL}, lemma 4.3] for the proof of the Reduction Theorem in this
article and for telling him about the paper \cite{Jo}. He is also grateful
to L.E. Jones for sending him a copy of the paper \cite{Jo}. The author
also thanks J.A. Hillman for telling him about some of his results related
to $PD_3$-groups.

\newpage
\bibliographystyle{plain}
\ifx\undefined\bysame
\newcommand{\bysame}{\leavevmode\hbox to3em{\hrulefill}\,}
\fi

\medskip

\end{document}